\documentclass[1p]{elsarticle}
\usepackage{graphicx}
\usepackage{booktabs, makecell, multirow, threeparttable}

\usepackage{caption}
\begin{document}

\begin{frontmatter}

\title{A hybrid heuristic algorithm for the resource-constrained
project scheduling  problem}
 \tnotetext[t1]{The study was carried
out within the framework of the state contract of the Sobolev
Institute of Mathematics (project FWNF-2022-0019)}
%
%\titlerunning{An improved GA for the RCPSP}
% If the paper title is too long for the running head, you can set
% an abbreviated paper title here
%

\author[1]{Evgenii N. Goncharov} \ead{gon@math.nsc.ru}

%\author{Evgenii N. Goncharov\inst{1,2}\orcidID{0000-0001-6843-8971} }
%
%\authorrunning{E. Goncharov}
% First names are abbreviated in the running head.
% If there are more than two authors, 'et al.' is used.
%

\address[1]{Sobolev Institute of Mathematics, Siberian Branch of the \\
Russian Academy of Sciences, Novosibirsk, 630090, Russia}

%\affiliation[1]{organization={Sobolev Institute of Mathematics},
%addressline={prosp. Akad. Koptyuga, 4}, postcode={630090},
%city={Novosibirsk}, country={Russia}}

%\institute{Sobolev Institute of Mathematics, Novosibirsk, prosp.
%Akad. Koptyuga, 4, Russia \\
%\email{gon@math.nsc.ru}\\
%\url{http://www.math.nsc.ru/}
% }
%
%\maketitle              % typeset the header of the contribution
%
\begin{abstract}
This study presents a hybrid metaheuristic for the
resource-constrained project scheduling problem (RCPSP), which
integrates a genetic algorithm (GA) and a neighborhood search
strategy (NS).
The RCPSP consists of a set of activities that follow
precedence relationship and consume resources.
The resources are renewable, and the amount of the resources is limited.
The objective of RCPSP is to find a schedule of the activities to minimize the
project makespan.
The algorithm uses two crossovers in the GA and
two neighborhoods in the NS, as well as a resource ranking
heuristic.
The computational results with instances from the PSPLIB
library validate the effectiveness of the proposed algorithm.
We have obtained some of the best average deviations of the solutions
from the critical path lower bound.
The best heuristic solutions have been updated for some instances
from PSPLIB.

\end{abstract}
\begin{keyword}
\texttt Project management  \sep resource-constrained project
scheduling problem \sep renewable resources \sep genetic algorithm
\sep neighborhood search \sep PCPLIB
\MSC[2020] 05C09 \sep 05C12 \sep 05C92
\end{keyword}

\end{frontmatter}

\section{Introduction}
The Resource-Constrained Project Scheduling Problem (RCPSP) may be
stated as follows: a project consists of a set of activities $N$
where each activity has to be processed without interruption.
A partial order on the set of activities is defined by a directed
acyclic graph.
The duration, the set, and the amounts of consumed
resources are assumed to be known for every activity.
The resource availability is assumed to be constant at each unit time interval
within the planning horizon of $\hat T$.
The resources outside the project horizon $\hat T$ considered to be unlimited.
All resources are renewable.
The objective is to schedule the activities of a project
to minimize the project makespan.

The RCPSP is a generalization of the Job Shop Scheduling problem.
This problem denoted as $m,1|cpm|C_{\max}$ using the classification scheme of
Herroelen et al. \cite{Herroelen1998a} or as $PS\ |\ prec\ |\ C_{\max}$
using the scheme of Brucker et al. \cite{Brucker1999}.
The RCPSP is NP-hard in the strong sense \cite{Blazewicz}.
It may be conceivable to use exact optimal methods only for projects of small size.
For the large-size instances, one needs heuristics to get the best solution
within a convenient response time, and heuristics remain the best
way to solve these problems efficiently.
It is worth noting that the RCPSP with cumulative resources can be
solved with polynomial complexity \cite{GimSev2003}.

The RCPSP is an important and challenging
problem that has been widely studied over the past
few decades.
Some of the most promising directions for developing heuristic methods
are based on genetic algorithms, local search methods,
bee and ant colony optimization, scatter search, and others.
We refer to the surveys provided by
Abdolshah (2014) \cite{Abdolshah2014},
Vanhoucke (2012) \cite{Vanhoucke2012},
Kolisch and Hartmann (2006) \cite{KolHar2006}.

Several new hybrid approaches have been put forward.
Rather than purely following the concept of one
single metaheuristic, these methods combine various
metaheuristic strategies, components,
and other optimization techniques.
A survey of hybrid metaheuristics for the
RCPSP provided by Pellerin at al. (2020) \cite{Pellerin2020}.
Researchers are also working on algorithms that integrate metaheuristics
and exact techniques to solve this problem \cite{VanhouckeCoelho2024}.
The hybridization of these two alternative
solution methods is known as a matheuristic.
Many extensions to the RCPSP have been proposed,
and an overview can be found in Hartmann and Briskorn (2010, 2021)
\cite{HarBris2010}, \cite{HarBris2021}.

The combination of GA algorithms with neighborhood search metaheuristics
is one of the most popular ways to hybridize metaheuristics.
We introduce a hybrid algorithm of this type.
It is based on the previously developed GA algorithm
\cite{Gon2017}, \cite{Gon2022a} and NS operator \cite{Gon2022}.
The GA algorithm uses crossovers in which the total resource
utilization of the current schedule is as high as possible.
We also emphasize the significance of unequal resources importance
and aim to maximize the use of the more valuable (scarce) resources.
As Coelho and Vanhoucke \cite{CoelhoVanhoucke2020} noted,
solving resource conflicts lies at the heart of the RCPSP,
and quite a number of studies have shown that the presence
of resources with limited availability is the main driver
of the hardness of the problem.
They modifies resources to find harder instances, whereas here the relative
importance of different types of resources is adjusted for better solving instances.
These concepts are related, our procedure should have an advantage over others
that do not consider resource importance.
We can also ''play'' with the weighting function of resource importance.
This expands the possibilities of the procedure, especially in cases
if all resources are approximately equally important.
A resource ranking is performed using a relaxed problem
with cumulative resources.

The proposed algorithm was tested on standard instance
datasets j60, j90, and j120 from the electronic library
PSPLIB \cite{PSPLIB}.
We have gotten the best average percent
deviations of the obtained solutions from the critical
path lower bound.
We have improved the best heuristic solutions for some instances
(4 instances from the dataset j90 and several dozens instances
from the dataset j120).
The results of numerical experiments are presented.

\section {Problem Setting}
The RCPSP problem can be defined as follows.
A project is taken as a directed acyclic graph $G=(N,A)$.
Denote by $ N =\{ 1, ..., n \} \cup \{ 0, n+1 \} $
the set of activities in the project where activities
$0$ and $n+1$ are dummy.
These dummy activities define the start and the
completion of the project, respectively.
The precedence relation on the set $N$ is defined with a set of pairs
$A = \{ (i, j) ~ | ~ i ~precedes  ~j \}$.
If $(i, j) \in A$, then activity $j$ cannot start before activity $i$
has been completed.
The set $A$ contains all pairs $(0, j)$ and $(j, n+1)$, $j=1, ..., n$.

Denote by $ K $ the set of renewable resources.
For each resource type $ k\in K $, there is
a constant availability $R_k \in Z^{+}$ throughout
the project horizon $ \hat T $.
Activity $j$ has deterministic duration $p_j\in Z^{+}$.
The profile of resource consumption is
assumed to be constant for every activity.
So, activity $j$ requires  $r_{jk} \geq 0$ units of resource
of type $k,\ k \in K$, at every time instant when it is processed.
We assume that $r_{jk} \leq R_k, \ j \in N, \ k \in K$.

Let's introduce the problem variables.
Denote by $s_j \geq 0$ the starting time of activity $j\in N$.
Since activities are executed without preemptions, the completion time of activity
$j$ is equal to $c_j = s_j + p_j$.
Define a schedule $S$ as an $(n+2)$-vector $(s_{0},...,s_{n+1})$.
The completion time $T(S)$ of the project corresponds to the moment
when the last activity ${n+1}$ is completed, i.e., $T(S) = c_{n+1}$.
Denote by $J(t) = \{ j \in N ~ | ~ s_j < t \leq c_j \}$
the set of activities which are
executed in the unit time interval $[t-1, t)$ under schedule $S$.
The problem is to find a feasible schedule $S=\{s_j \}$
respecting the resource and precedence constraints so that
the completion time of the project is minimized.
It can be formalized as follows: minimize the makespan of the project
\begin{equation}
\label{e1} T(S)=\max_{j\in N} (s_{j}+p_{j})
\end{equation}
under constraints
\begin{equation}
\label{e2} s_{i}+p_{i}\leq s_{j}, ~~\ \forall (i, j) \in A;
\end{equation}
\begin{equation}
\label{e3} \sum\limits_{j \in J(t)}r_{jk} \leq R_k ,~ k \in K,
~ t=1,...,\hat T;
\end{equation}
\begin{equation}
\label{e4} s_j\in {\bf Z}^+,  \quad j\in N.
\end {equation}
The objective function (\ref{e1}) minimizes the makespan of the project.
Constraints (\ref{e2}) enforce the precedence constraints between activities.
Relation (\ref{e3}) corresponds to the resource constraints,
and condition (\ref{e4}) defines the decision variables.

\section {Resource ranking}
Coelho and Vanhoucke \cite{CoelhoVanhoucke2020} analyzed the hardness
of RCPSP from an experimental point-of-view
and proposed an approach to transform easy RCPSP instances
into very hard ones.
These instances should be as small as possible in terms
of the number of activities and resources.
One of the steps in this procedure is to reduce the number
of resources needed to search for bottleneck resources that
really define the instance's hardness.

We conducted numerical experiments on standard PSPLIB instances.
Of course, we do not remove resources from these instances.
However, we try to predict less significant resources and reduce their
influence on the solution, while increasing the influence of more
significant ones.
In our opinion, this approach echoes the approach
of Coelho and Vanhoucke.

We use the concept of greater or lesser degrees
of ''scarcity'' for resources.
Resources are ranked according to their degree of scarcity.
The ranking of resources can be obtained by solving a relaxed
problem with cumulative resources.

Let us consider a relaxed problem with cumulative resources:
\begin{equation}
\label{e5} T(S)=\max_{j\in N} (s_{j}+p_{j})
\end{equation}
under constraints
\begin{equation}
\label{e6} s_{i}+p_{i}\leq s_{j}, ~~\ \forall (i, j) \in A;
\end{equation}
\begin{equation}
\label{e7} \sum_{t'=1}^{t} \sum_{j\in J(t')} r_{jk} \leq
\sum_{t'=1}^{t}R_{k}, ~ k \in K, ~ t=1,...,\hat T;
\end{equation}
\begin{equation}
\label{e8} s_j\in {\bf Z}^+,  \quad j\in N.
\end {equation}
The relation (\ref{e7}) ensures the fulfillment of resource constraints
with cumulative resources, and (\ref{e5})--(\ref{e8})
is a scheduling problem with cumulative resources.

The problem (\ref{e5})--(\ref{e8}) can be resolved with
a known fast approximate method.
This algorithm is asymptotically exact with the deviation
tending to zero as the problem dimension increases \cite{Gim}.
An exact algorithm is also known for this relaxed problem,
where the activity duration is an integer \cite{GimSev2003}.
We use a fast approximate algorithm.
The complexity of this algorithm is $O(n^2 \log n)$.
As a result of using this algorithm, we obtain the remaining
cumulative resources, as well as the approximate solution for the relaxed problem.
We will consider resources with less unused balance to be more important (more scarce).
We hypothesize that the ranking of resources in
the relaxed problem is also conserved in the original RCPSP problem.
Let's renumber the resources in order of their ranks
and assign weights $ w_{k}, k \in K $ to resources
according to their ranks.
These weights can be defined in different ways (see Section 8).

\section {Genetic Algorithm}
We use the previously developed genetic algorithm for the RCPSP
problem \cite{Gon2017,Gon2022a}.
The process employs two crossovers.
The algorithm uses a heuristic rule to select promising segments (genes) from
the parent chromosomes, with the goal of incorporating these segments
into the offspring.
We call these genes as dense genes and their total resource
utilization during the current schedule is as high as possible.
Let's briefly describe the GA algorithm
(we refer to \cite{Gon2017,Gon2022a} for more details).

Represent a feasible solution as a list of activities
$ L = (j_0, ..., j_ {n+1}) $.
We use two algorithms for constructing schedules
of activities with serial (S-SGS) and parallel (P-SGS) decoders.
The serial decoding procedure \cite{KolHar1999} calculates an active
schedule $S(L)$ for an arbitrary list $L$.
A schedule is called active if none of the activities can be
started earlier without delaying some other activity.
It is known that there is an optimal schedule among active schedules.
The parallel decoder (P-SGS) processes schedule time in sequence,
identifying all the activities that can be initiated prior to the scheduled time.

Each schedule for the {\it initial population} $\Gamma $
is constructed as follows.
We generate a random feasible list of activities $L$
and then construct a schedule $S$ with the parallel decoder.
Further, we apply a local improvement procedure FBI
(the forward-backward improvement procedure) to the resulting schedule.
The best schedule found through this process is then incorporated
into the population $\Gamma$.
 These steps are repeated until $\Gamma $ contains the necessary
number of schedules.

{\bf The set of parent schedules.}
We look through the entire population $\Gamma $ in
nondecreasing order of schedule length and add each
schedule to the {\it set of parents} $\Gamma^\prime$ with a
certain probability.
If we have looked through the entire population and have not collected the
necessary number of schedules $ParentsSize$ (an algorithm parameter),
we choose the best among those which have not yet been added.
A pair of parent schedules will be randomly selected
for a crossover from the set of parent schedules $\Gamma^\prime$.
This stra\-tegy for choosing two parent schedules has been shown to be effective
in problems of high dimension (see for instance \cite{ProonJin2011}).

We transmit the segment of the parent chromosome (gene) to the offspring
chromosome, and we aim to ensure that this gene is in high resource
utilization of the current schedule.
At every scheduled time, we discover a surplus of unused resources.
When the residue is less than a predefined admissible residue,
we determine the set of activities accomplished during this time interval (a dense gene).
In general, we can operate with the weighted residues of unused resources.

Let $\Theta $ be the set of dense genes.
Comparing gene weights allow us to prioritize those that utilize
high-priority resources more rationally, leading to a decrease
in surplus unused resources.

\emph{ $DenseActivities(S, R) \rightarrow \Theta$. }
\begin{enumerate}
    \item Set $t := 0$, $\Theta = \emptyset$.
    \item  While $ t < T(S) $:
    \begin{enumerate}
        \item  find a weight $v_t$ of the set $J(t)$:~
        $v_t=\sum_{k\in K} \left ( R_{k}-\sum_{j\in J(t)}r_{jk} \right )w_{k}/R_{k}  $;
        \item if $ v_t < R $,  then update the list of dense genes $\Theta = \Theta \cup J(t)$;
        \item  $ t := t + 1 $.
    \end{enumerate}
\end{enumerate}
It may happen that some dense genes overlap.
In this case, we will keep the gene that produces the smaller
weighted surplus of unused resources $v_t$.

\section {Crossovers}
Once a parent chromosomes have been selected, a crossover operator will
be applied on the chromosomes to construct offspring chromosomes.
Two crossovers have been constructed.

{\bf Crossover procedure A.}
We consider the first dense genes in each parent schedule,
and the leading gene should be the gene with the lowest weight.
All the initial activities from the parent chromosome, including the
leading gene, are incorporated into the offspring chromosome.
We then remove this leading gene from the list of dense genes.
Let's look at the first dense genes in the chromosomes that haven't had
any activities added to a daughter chromosome.
We choose the best of the two new genes and add the activities from
the parental chromosome that contained them, excluding the activities
that have already been added to the daughter chromosome.
These steps are repeated until we have considered all dense genes.
If it's necessary, we add the remaining activities to
the daughter chromosome in the same order as they are located
in the chromosome with minimal duration.
The serial decoder is used to construct the daughter schedule
from the daughter chromosome.

{\bf Crossover procedure B.}
 A block of activity $j$ in an active schedule $S$ is a set
of activities that overlap $j$ in this feasible schedule,
start immediately after $j$, or finish immediately before it.
Graph $G_S$ is the directed graph with the vertex set $V = N$ and
the arc set $E = \{ (i, j) ~ | ~ c_i = s_j, (i, j) \in A \}$.
The outgoing network of activity $j$ for schedule $S$ is the maximal
(by inclusion) connected subgraph of graph $G_S$,
with the only source being the vertex that corresponds to activity $j$.
We select the dense genes with the lowest weight criteria
from the both parent chromosomes.
Identify these activities in the second chromosome and
find the outgoing (or incoming) network for each of them.
Then find a segment in the list of activities
between the leftmost and rightmost activities in the block
and outgoing (incoming) networks.

{\it The mutation operator} involves ''mixing'' activities
from the list corresponding to the schedule and generates a
schedule from the resulting combined list of activities.
The algorithm applies the mutation operator on chromosomes to avoid
being trapped in local optimum solutions.
The mutation is carried out in two phases.
The first step involves selecting two random activities in the chromosome
and swapping them if they do not violate the priority conditions.
In the second step, we relocate a random activity to a different
location without violating the priority conditions.
This process is iterated a specified number of times.

We pick a given number (a parameter of the algorithm)
of the top chromosomes from the set of offspring chromosomes
$\Gamma^{\prime\prime}$.
They are included in {\it the next generation} $\Gamma $.
The same number of the worst chromosomes are removed from $\Gamma $.

\section {Neighborhood search}
Let's briefly describe the NS operator (for more details,
we refer to \cite{Gon2022}).
For a given feasible schedule $S = (s_{0},...,s_{n+1})$ and a
core activity $j=1, ..., n$, the NS operator reschedules a set of activities
$A^s_j$, while maintaining the start times of all other activities.
Let $P$ be a predetermined number of activities to be rescheduled.
A value of $P$ affects the computation time to
obtain a neighborhood solution by rescheduling.
The smaller value of $P$ usually means fewer activities that need
to be rescheduled and less time to obtain a new schedule.
The following Block selection
method is used to create $A^s_j$ \cite{Palpant2004}.\\
\emph{ $CreateBlock(j,S) \rightarrow A^s_j$. }
\begin{enumerate}
    \item $A^s_j= {j};\ b = 0;$
 create a random order for all activities in $ A/\{j\}$.
 Let $i$ be the first activity in the order.
    \item If  $s_j - p_i - b \leq s_i \leq s_j + p_j + b,\ $ then $A^s_j= A^s_j\bigcup \{i\}.$
    \item If $|A^s_j| = P$, go to Step 6.
    \item If $i$ is the last activity among the ones not belonging
    to $A^s_j$ based on the order defined in Step 1, $b = b+1$.
    \item Let $i$ be the next activity among the ones not belonging to $A^s_j$
based on the order defined in Step 1. Go to Step 2.
    \item END.
\end{enumerate}
The Block selection method basically selects a set $P$ of activities
that overlap, or are close to an activity $j$ in a given feasible schedule.

The NS method uses knowledge gained from the previously considered solutions.
This knowledge is stored in a {\it tabu list}, which serves
as a memory to avoid cyclicity.
Each entry in the taboo list records the attributes
of the last visited solutions.
The length of a taboo list can vary based on the details of the NS operator's process.
As the taboo status of an arbitrary solution $L$,
we consider the sum of the start times:
$$
TS(L) = \sum_{j=1}^{n} s_{j}
$$
for the schedule $S(L)$.

The neighborhood search algorithm uses two types
of neighborhoods \cite {Gon2022}.
The first neighborhood $N_A(S)$ is a modification of the scheme
proposed in \cite{ProonJin2011}.
For a given feasible schedule $S = \{s_0,s_1,...,s_n,s_{n+1}\}$ and a
core activity $j \in A$, we define the activity block $A^s_j$.
The NS operator reschedules the set of activities $A^s_j$,
keeping the start times of the remaining activities intact.
The rescheduling sub-problem is formed by the following steps.
We fix the start time of all activities that are not part of the set
$A^s_j$ and release the resources used by all
activities from $A^s_j$  for each time period $t$.
The available amount of resource $k$ for activities
in $A^s_j$  in period $t$ is equal to $R_k$ minus
the resources used by all the activities in $A^s_j$ in period $t$.
Then we derive the earliest start time (EST) and latest finish time (LFT)
for each activity $i\in A^s_j$  as
$EST_i = max\{s_l + p_l,\ \forall\ l\not\in A^s_j \ and \ (l,i) \in A\}$ and
$LFT_i = max\{s_l,\ \forall\ l\not\in A^s_j \ and \ (l,i) \in A\}$.
Interval $(EST_i ,LFT_i )$ establishes the time window for activity $i$
that can be rescheduled to ensure that a new schedule remains
feasible for all remaining activities.
The rescheduling problem is to reschedule all activities in $A^s_j$
so as to minimize their makespan, while meeting the resource restrictions
of each period and the time window constraints defined by
$(EST_i ,LFT_i )$.

At each iteration, a new random vector is generated for the
activities of $A^s_j$  as a priority list.
We order the activities in $A^s_j$ by decreasing their weights $v_j$.
The vector is created iteratively by randomly selecting the next activity
from the ordered list among all unselected activities that preceded
the ones in $A^s_j$.
Each activity in the priority list is moved to the earliest
(latest) start time that is precedence and resource feasible, and satisfies
the time window $(EST_i ,LFT_i )$.
Following the rescheduling of all activities of $A^s_j$,
activities outside of set $A^s_j$ are incorporated to construct
a complete feasible solution.
All activities in $A$ undergo a global left shift to minimize
the makespan as much as possible.
A new schedule is compared against the previous solution prior
to applying the NS operator.
If the makespan is improved, the resulting schedule replace
the previous schedule.
If no improvement occurs prior to reaching the predefined limit $\lambda_{NS}$
of iterations, the S-SGS is then applied to the schedule with a fresh
random priority list at the next iteration.

The neighborhood $N_B(S)$ is defined as follows.
We determine the block of activity $A^s_j$ for a given list of activities
$L$ (and a correspondent active schedule $S = \{s_0,...,s_{n+1}\}$),
where $j$ is a core activity.
If the block contains at least one predecessor of activity $j$, then we set
$A^s_j$ to be empty.
The list $L$ is represented in the form of three
successive lists $L = A^1, A^s_j, A^2$.
For each activity $j \in A$, an element $L'$ of the
neighborhood $N_B(S)$ is constructed using a non-empty block $A^s_j$.
The list $L'$ is generated from the list $L$ by following these steps.
For set $A^1$, the start times of all activities are adjusted,
and the resources, occupied by all activities in set $A^1$, are
relinquished at each time period $t$.
A partial schedule for the activities in set $A^1$ is generated
by employing the serial decoding procedure.
We then expand the partial schedule by incorporating activities
from the set $A^s_j$, utilizing a parallel decoding technique.
According to the procedure, for each schedule time $t$ we have
the corresponding eligible set $E_t$, i.e. a set of activities
that can be started at time $t$ without violating any constraints.
There are exponentially number of possibilities to select a subset
of activities from the eligible set to include into the schedule.
We solve the multi-dimensional knapsack problem with an objective
function that maximizes the weighted resource utilization ratio:
\begin{equation}
\label{e9} max \ \sum_{j\in E_t} x_j \sum_{k\in K} \frac{w_{k}r_{jk}}{R_{k}},
\end{equation}
\begin{equation}
\label{e10} \sum_{j\in E_t} r_{jk}x_j \leq R_{k} - \sum_{j\in J(t)} r_{jk} \, \quad k\in K,
\end{equation}
\begin{equation}
\label{e11} x_j\in \{0,1\},  \quad j\in N.
\end {equation}
The right-hand side of the restriction represents the remaining capacity
of resource type $k$ at time $t$.
The problem is resolved using Greedy Randomized Adaptive Search Procedures (GRASP).

Finally, we construct the list $L'$ as follows.
The activities from $A^s_j$ are then added to the list $L'$ in ascending
order of their start times within the partial schedule.
The remaining activities are ordered similarly in the list $L'$ as in the list $L$.
The schedule $S(L')$ is referred to as a neighboring sample for the schedule $S$.
The set comprising all adjacent samples is referred to as the neighborhood
of the schedule $S$ and is denoted by $N_B(S)$.

\section {The GANS Algorithm}
We propose a genetic algorithm for solving the RCPSP
by integrating a neighborhood search approach.
The proposed GANS algorithm seeks to preserve the randomness of the
Genetic Algorithm search and enhance solution quality
through a neighborhood search conducted at specific points
during the GA iterations.

Numerical experiments were conducted using instances
from the PSPLIB library.
Separate analyses for the GA and NS algorithms had previously
revealed varying degrees of efficiency across different groups
of instances from PSPLIB.
Let us find deviation of current upper bound from its critical path value
for each instance from a group in the j120 and compute the average deviation
for each group.
Table 1 presents the top of the sorted list of instance groups
in non-decreasing order of their average deviations.
For each group, the average deviations of solutions, obtained
with the GA and NS algorithms, are shown.
The best performing results are highlighted.
\begin{table}[!ht]
\caption{Average percent deviations from the critical path value
  for the some groups of j120.}\label{tab4}
    \centering
    \begin{tabular}{|l|l|l|}
    \hline
        Instance Set & GA operator \cite{Gon2017} & NS operator \cite {Gon2022}\\ \hline
        J120\_56 & 156,83 & {\bf 154,51}  \\ \hline
        J120\_16 & 145,39 & {\bf 143,90}  \\ \hline
        J120\_36 & 131,14 & {\bf 129,37}  \\ \hline
        J120\_51 & 125,85 & {\bf 123,64}  \\ \hline
        J120\_31 & 112,33 & {\bf 111,15}  \\ \hline
        J120\_11 & 108,15 & {\bf 107,16}  \\ \hline
        J120\_6  & {\bf 76,97} & 77,64  \\ \hline
        J120\_46 & {\bf 73,09} & 73,46  \\ \hline
        J120\_26 & 69,05 & {\bf 68,66}  \\ \hline
        J120\_37 & {\bf 66,72} & 66,94  \\ \hline
        J120\_52 & 55,95 & {\bf 55,40}  \\ \hline
        J120\_57 & {\bf 66,35} & 66,41  \\ \hline
        J120\_17 & {\bf 60,64} & 61,75  \\ \hline
    \end{tabular}
\end{table}
The performance of GA is equal to or better than NS operator for all
remaining groups of instances not listed in the table.

A full factorial experimental design of instances is defined
by combining three parameters: network complexity (NC),
resource factor (RF), and resource strength (RS).
The NC defines the average number of precedence relations per activity.
The resource factor RF reflects the average portion of resources requested per activity.
The RS measures scarcity of the resources.
A zero RS factor indicates the lowest amount of each resource type needed
to complete all activities, whereas an RS value of one represents the required
amount of each resource type allocated according to the early start time schedule.
It is known \cite {Weglarz1999} that values of the parameters
$RF = 4$, $RS = 0.1$ or $RS = 0.2$ match hard enough series.
As shown in Table 1, a neighborhood search is more preferable
on challenging instances with a small $RS$ value.

This observation has led to the idea of dividing the
PSPLIB's instances into three subsets.
Instances from these subsets will be processed using the same hybrid algorithm,
but with a separate set of problem parameters.
This will allow us to efficiently coordinate the problem-solving process
within each of the three subsets of instances.
Let us establish the threshold values $\sigma_{1}$ and $\sigma_{2}$
to define subsets, where $0<\sigma_{1}<\sigma_{2}$.
The process starts with the formation of the initial population
by the Genetic Algorithm.
Then we find the best chromosome in the initial population and calculate
the value of the relative deviation of its makespan from the critical path length
 $\sigma = (UB - Tkr)/Tkr$.
If $\sigma < \sigma_{1}$, then the instance is placed in the first subset.
If $\sigma_{1} \leq \sigma \leq \sigma_{2}$, then the instance is included
in the second subset.
If $\sigma \geq \sigma_{2} $, then the instance is added to the third subset.
The primary distinction in applying the algorithm to each subset of instances
lies in the varying proportions of GA and NS operators application.
The first subset mainly employs the GA operator, whereas the third
subset primarily utilizes the NS operator.
The operators in the second subset have a proportion that is roughly equivalent.

The GA starts running and will persist until the count of its iterations
reaches a predetermined value without modifying the existing record (best solution).
In this instance, the NS operator initiates.
The number of schedules generated by a single NS operator is a parameter of the problem.

The initial solution in the NS operator can be chosen using any suitable
method that is available.
The choice of the algorithm for the initial solution is not critical
for the local search methods.
Gagnon at al. \cite{Gagnon2004} noted that there is some dilemma concerning
the choice of the initial solution used by a NS method adaptation.
Starting with a very good solution doesn't let enough space to find
a significant improvement.
On the other side, it may take a long computation time to improve a bad starting
solution.
An initial solution will be chosen from the set of chromosomes in the population at random
using a specified probability -- a predefined parameter, and the selection process
will begin with the best chromosome.

Denote by $\lambda $ the maximum number of schedules generated.
The number of schedules generated will increase every time the objective function is
computed, and the algorithm stops when $\lambda $ is exceeded.
The general scheme of the algorithm is outlined below.

\begin{enumerate}

    \item  Solve the relaxed problem (\ref{e5})--(\ref{e8}).
       Rank resources, renumber them in a descending order of their scarcity,
       and assign resource weights $ w_{k}, k \in K $.
    \item  Construct the initial population $\Gamma$ and store the best of its chromosomes as $S^\ast$ along with its makespan $T^*$.
    \item  Calculate  $\sigma$ by setting $UB$ equal to $T^*$. Then we
    decide to categorize the instance within one of the three subsets of the set of instances PSPLIB.
    Assign parameters of the problem based on the selected option.
    \item  Initiate a set of offspring schedules $\Gamma^{\prime\prime} := \emptyset$.
    \item While the number of schedules generated does not exceed $\lambda $, do:
    \begin{enumerate}
        \item construct the set of parent chromosomes $\Gamma^\prime $,
        \item until the necessary number of offspring chromosomes has been generated, do:
        \begin{enumerate}
            \item choose two parent chromosomes $S_1$ and $S_2 \in \Gamma^\prime$ at random,
            \item choose crossover $Crossing$ with equal probability from $CrossingA$ and $CrossingB$,
            \item cross the parent chromosomes $S_1$ and $S_2$ :  $S^\prime := Crossing(S_1, S_2) $,
            \item apply to $S^\prime$ sequentially the operations of mutation and local improvement FBI,
            \item if the mutation has made $S^\prime$ worse, it is canceled,
            \item if $T(S^\prime) < T(S^\ast)$,  then update the record $ S^\ast := S^\prime$,
            \item update the set of offspring chromosomes
              $\Gamma^{\prime\prime} = \Gamma^{\prime\prime} \cup  S^\prime  $,
        \end{enumerate}
        \item create a new population for the next generation,
        \item if $S^\ast$ has not been updated for a given number of steps, then
         \begin{enumerate}
         \item replace a given number of the worst schedules in the population with new chromosomes;
         \item  Execute a particular number of steps (an algorithm parameter) in the NS operator:
        \begin{enumerate}
            \item choose the initial chromosomes $S$, and set  $T: = T(S)$.
Tabu list $TL$ is set empty,
            \item Until the stopping criterion is satisfied, do:
        \begin{itemize}
            \item choose the type of neighborhood ($N_A(S)$ or $N_B(S)$) equally probable,
            \item Find the neighbor sample $S'$, not prohibited by the tabu list $TL$,
            \item If $\ T(S') < T$, then assume $T: = T(S'),\ S: = S'$,
            \item Update the tabu list $TL$, and set $S: = S'$.
        \end{itemize}
            \item If $\ T < T^*$, then assume $T^*: = T(S'),\ S^*: = S'$,
         \end{enumerate}
        \end{enumerate}

    \end{enumerate}

\end{enumerate}
The chromosome $S^\ast$ is the result of the algorithm.

\section {Search intensification and parameters}
The parameters of the problem exert a substantial influence on determining its solution.
We adjust some parameters during the problem-solving process, introducing
self-tuning elements into the algorithm.

Computational experiments indicate that $\sigma_{1}$ should be set
to $0.2$ and $\sigma_{2}$ to $0.6$.
We use a population size $\Gamma $ ranging from 40 to 150 chromosomes.
The population size is decreased for instances produced with a low NC value,
and for a high NC value, the population size is increased.
We assume an average population size of 60-80 chromosomes.
If we expend considerable effort to replenish the population with new unique chromosomes,
we will either decrease the population size or waive the requirement for chromosome
uniqueness.
Otherwise, we may be able to increase the population size.

The resource weights $ w $ are chosen with equal probability from the following set of options:
$ (1,0.8,0.6,0.4)$, $ (1,0.9,0.8,0.7)$,  $ (1,1,1,1)$,
and $ w_{k} = 2 - \frac{\hat{R}_k}{\hat{R}_4} ,\ k \in K$.
In the last case, the reduction in the value of $ w_{k}$ is directly proportional
to the remaining unused resources $\hat{R}_k,\ k \in K$ in the relaxed problem.
This version of weighting coefficients was also applied in a probabilistic implementation
for the unlimited scenario.
The parameter $R$ in the operator $DenseActivities(S, R)$ plays a key role in
determining the quality of the dense gene that is being identified.
Relatively small $R$ enable the discovery a high-quality dense gene.
It's possible that we may not discover these genes or there are only a limited number of them.
We aim to enhance the quality of the dense gene, which is achieved by reducing $R$,
but in this case, we are compelled to increase it.
Conversely, if the cardinality of the dense genes set is sufficiently great,
then the algorithm can decrease the parameter $R$.
We assume the initial value as $R=0.75$.

A higher probability of a chromosome being included in the set of parent chromosomes
results in better overall quality, however it may not be sufficient to produce
offspring with sufficient genetic diversity.
Under certain conditions, the algorithm has the capability to alter this parameter,
either by reducing or enhancing it.
We assume the initial value as $probabilityParentSelection=0.25$.

If the current record $S$ (the best chromosome) in the operator NS remains unchanged
over a specified number of iterations, we modify the cardinality of the
set $A^s_j$ within the procedure $CreateBlock(j,S)$.
This cardinality can be either increased or decreased based on the
statistics of the previous iterations of the NS algorithm.
If a proportion of non-empty neighboring solutions in a specified number
of previous iterations meets or exceeds a predetermined threshold,
we then increase the parameter $P$ by one.
A high proportion of empty neighbors implies that a predecessors of activities
being frequently included in the set $A^s_j$.
This means that the cardinality of this set is excessively high,
we decrease the parameter $P$ by one in this case.
The intention behind such an adjustment is to enhance
the size of the set $ A^s_j$, with the goal of yielding new
solutions.
This increase can only proceed up to a certain limit: when the set
often includes a predecessors and is forced to become empty,
it is crucial to decrease the parameter $P$, thereby re-entering
the domain of non-empty sets $A^s_j$.
If the parameter $P$ has undergone a fixed number of
modifications and the value of $T$ remains unchanged, then reset
$P$ to 1 and reassign the resource weights $ w_{k}, k \in K $.

\section { Numerical experiments.}
The GANS algorithm was implemented in C++ using Visual Studio and executed on a computer
featuring a 3.4 GHz CPU and 16 GB of RAM with the Windows 10 operating system.
In evaluating the performance of the proposed algorithm, we employ the standard set
of instances provided by Kolisch and Sprecher \cite{PSPLIB}.
The project scheduling library PSPLIB contains these instances along
with their best-known values.
The instances are downloadable at http://www.om-db.wi.tum.de/psplib/.

The optimal solutions for the instances from datasets j60, j90, and j120 are unknown.
The quality of the solutions is gauged by the average percent deviation (APD)
of the received solutions from the lower bounds calculated by the critical path method.
The stopping criterion is determined by the maximum number $\lambda $ of schedules
generated, as per Kolisch and Hartmann \cite{KolHar2006}.
The selected limits on the number of generated schedules are 50000, 500000, and unlimited.

Comparisons between the performance of the GANS algorithm and
previous results from experimental evaluations of competitive
heuristics for the datasets j60, j90, and j120 are shown in Tables
\ref{tab1}--\ref{tab3}. The presented results confirm the high
effectiveness of the proposed algorithm through a clear analysis.

\begin{table}
\footnotesize \caption{Average percent deviations from the critical path
 lower bound for the dataset j60.}\label{tab1}
\resizebox{\textwidth}{!}{%
\begin{tabular}{p{100pt} p{124pt}  c c c}
\hline & & &
\multicolumn{2}{p{110pt}}{\ \ \ \ \ \ \  APD, \%}  \\
Algorithm & Reference & $\lambda =50000$ & $\lambda =500000$ &$ unlimited $\\
\hline
GANS               & this paper (2024)                                 & {\bf 10,48} & {\bf 10,37}  & {\bf 10,37} \\
GA                 & Goncharov \cite{Gon2022a} (2022)                  & 10,50 & 10,40  &  \\
GA                 & Goncharov, Leonov\cite{Gon2017} (2017)            & 10,52 & 10,42  &   \\
GANS               & Proon, Jin \cite{ProonJin2011} (2011)             & 10,52 & -- & 10,52 \\
TS,VNS             & Goncharov \cite{Gon2022} (2022)                   & 10,55 &  10,44 & \\
Sequential(SS(FBI))& Berthaut et al. \cite{Berthaut2018} (2018)        & 10,58 & 10,45 & \\
GH + SS(LS)        & Paraskevopoulos et al. \cite{Paraskevopoulos2012} (2012) & 10,54 & 10,46 & \\
AI(FBI)            & Mobini, at al. \cite{Mobini2011} (2011)           & 10,55 & -- & \\
TS + SS(FBI)       & Mobini, at al. \cite{Mobini2009} (2009)           & 10,57 & -- & \\
GA(FBI)            & Wang et al. \cite{WangLi2010} (2010)              & 10,57 & -- & \\
GA(FBI)            & Goncalves \cite{Goncalves2011} (2011)             & 10,57 & 10,49 & \\
EA(GA(LS)+DEA(LS)) & Elsayed et al. \cite{Elsayed2017} (2017)          & 10,58 & -- & \\
PSO(LS)            & Czogalla and Fink \cite{Czogalla2009} (2009)      & 10,62 &  -- & \\
GA                 & Lim et al. \cite{Lim2013} (2013)                  & 10,63 &  10,51 & \\
Parallel(MA(LS))   & Chen, at al. \cite{ChenLiu2014} (2014)            & 10,63 & -- & \\
PL(LS)             & Zheng and Wang \cite{ZhengWang2015} (2015)        & 10,64 & -- & \\
GA(FBI)            & Zamani \cite{Zamani2013} (2013)                   & 10,65 & -- & \\
SFL(LS)            & Fang and Wang \cite{FangWang2012} (2012)          & 10,66 & -- & \\
GA(FBI)            & Ismail and Barghash \cite{Ismail2012} (2012)      & 10,66 & -- & \\
ACOSS              & Wang Chen, at al. \cite{WangChen2010} (2010)      & 10,67 & -- & \\
GAPS               & Mendes, at al. \cite{Mendes2009} (2009)           & 10,67 & 10,67 & \\
GA                 & Debels, Vanhoucke \cite{DebelsVanhoucke2007} (2007) & 10,68 & -- & \\
Specialist(PSO(LS))& Koulinas et al. \cite{Koulinas2014} (2014)        & 10,68 & -- & \\
GA(LS)             & Carlier et al. \cite{Carlier2009} (2009)          & 10,70 & -- & \\
Decomposition      & Palpant et al. \cite{Palpant2004} (2004)          & -- & -- & 10,81 \\
Population-based   & Valls et al. \cite{Valls2004} (2004)              & -- & -- & 10,89 \\
\hline
\end{tabular}
              }
\end{table}

\begin{table}
\footnotesize
\caption{Average percent deviations from the critical path
 lower bound for the dataset j90.}\label{tab2}
\resizebox{\textwidth}{!}{%
\begin{tabular}{p{94pt} p{130pt}  c c c}
\hline & & &
\multicolumn{2}{p{110pt}}{\ \ \ \ \ \ \  APD, \%}  \\
Algorithm & Reference & $\lambda =50000$ & $\lambda =500000$ &$ unlimited $\\
\hline
GANS               & this paper                                     & {\bf 9,80} &  {\bf 9,53} & {\bf 9,48} \\
GA                 & Goncharov \cite{Gon2022a} (2022)               & 9,92 & 9,61 &  \\
GA                 & Debels, Vanhoucke \cite{DebelsVanhoucke2007} (2017) & 9,90 & -- & \\
Sequential(SS(FBI))& Berthaut et al. \cite{Berthaut2018} (2018)     & 9,96       & 9,74 & \\
TS,VNS             & Goncharov \cite{Gon2022} (2022)                & 9,98       &  9,78 & \\
Sequential(SS)     & Ranjbar and Kianfar \cite{Ranjbar2009} (2009)  & 10,04      & -- & \\
SS(EM + FBI)       & Debels, et al. \cite{Debels2006} (2006)        & 10,09      & 9,80 & \\
PL(LS)             & Jedrzejowicz, Ratajczak \cite{Jedrzejowicz2006} (2006) & 11,60   & -- & \\
TS                 & Ying et al. \cite{YingLinLee2009} (2009)       & 12,15      &  -- & \\
\hline
\end{tabular}
              }
\end{table}

\begin{table}
\footnotesize
\caption{Average percent deviations from the critical path
 lower bound for the dataset j120.}\label{tab3}
\resizebox{\textwidth}{!}{%
\begin{tabular}{p{100pt} p{124pt}  c c c}
\hline & & &
\multicolumn{2}{p{110pt}}{\ \ \ \ \ \ \  APD, \%}  \\
Algorithm & Reference & $\lambda =50000$ & $\lambda =500000$ &$ unlimited $\\
\hline
GANS                 & this paper                             & {\bf 30,42} &  {\bf 29,37}  & {\bf 29,19} \\
GA                   & Goncharov \cite{Gon2022a} (2022)       & 30,46 & 29,63 &  \\
Specialist GA        & Goncharov, Leonov\cite{Gon2017} (2017) & 30,50 &  29,74 & \\
TS,VNS               & Goncharov \cite {Gon2022} (2022) & 30,56 &  29,88 & \\
GA                   & Lim et al. \cite{Lim2013} (2013)       & 30,66 &  29,91 & \\
biased random-key GA & Goncalves \cite{Goncalves2011} (2011)  & 32,76 & 30,08 & \\
GANS                 & Proon, Jin \cite{ProonJin2011} (2011)  & 30,45 & 30,78 & 30,78 \\
ACOSS                & Wang Chen, at al. \cite{WangChen2010} (2010) & 30,56 & -- & \\
DBGA                 & Debels, Vanhoucke \cite{DebelsVanhoucke2007} (2007) & 30,69 & -- & \\
GH + SS(LS)          & Paraskevopoulos et al. \cite{Paraskevopoulos2012} (2012) & 30,78 & 30,39 & \\
GA                   & Debels, Vanhoucke \cite{DebelsVanhoucke2007} (2007) & 30,82 & -- & \\
PL(LS)               & Zheng and Wang \cite{ZhengWang2015} (2015) & 31,02 & -- & \\
SFL(LS)              & Fang and Wang \cite{FangWang2012} (2012) & 31,11 & -- & \\
Sequential(SS(FBI))  & Berthaut et al. \cite{Berthaut2018} (2018) & 31,16 & 30,39 & \\
EA(GA(LS)+DEA(LS))   & Elsayed et al. \cite{Elsayed2017} (2017) & 31,22 & -- & \\
Specialist(PSO(LS))  & Koulinas et al. \cite{Koulinas2014} (2014) & 31,23 & -- & \\
GA - Hybrid, FBI     & Valls, at al. \cite{Valls2008} (2008)  & 31,24 & 30,95 & 30,95 \\
GA(FBI)              & Wang et al. \cite{WangLi2010} (2010)  & 31,28 & -- & \\
GA(FBI)              & Zamani \cite{Zamani2013} (2013)       & 31,30 & -- & \\
Enhanced SS          & Mobini, at al. \cite{Mobini2009} (2009) & 31,37 & -- & \\
GA(LS)               & Alcaraz and Maroto \cite{Alcaraz2006} (2006) & 31,38 & -- & \\
GA(LS)               & Carlier et al. \cite{Carlier2009} (2009) & 31,40 & -- & \\
Scatter search - FBI & Debels, et al. \cite{Debels2006} (2006) & 31,57 & 30,48 & \\
GAPS                 & Mendes, at al. \cite{Mendes2009} (2009) & 31,44 & 31,20 & \\
GA, FBI              & Valls, et al. \cite{Valls2005} (2005) & 31,58 & -- &\\
Decomposition        & Palpant et al. \cite{Palpant2004} (2004)          & -- & -- & 31,58 \\
\hline
\end{tabular}
              }
\end{table}

For the dataset j60, the best known solutions (at the time of writing)
were achieved for all instances except one (j609\_10).
Note that for j60 the heuristic solutions have remained unchanged since 2008.
The best APD values were achieved for datasets j90 and j120 at every lambda value.
We have improved the  existing best heuristic solutions for some instances
(4 instances from the dataset j90 and several dozens instances
from the dataset j120, see http://www.om-db.wi.tum.de/psplib/).

The links to text files containing the achieved makespans for all instances
can be found in Table \ref{tab5}.
\begin{table}[!ht]
\caption{Hyperlinks to the files with makespans.}\label{tab5}
    \centering
    \begin{tabular}{|l|l|}
    \hline
        dataset & link to file \\ \hline
        j60 &  old.math.nsc.ru/LBRT/k4/j60-hrs-Goncharov.txt \\ \hline
        j90 &  old.math.nsc.ru/LBRT/k4/j90-hrs-Goncharov.txt \\ \hline
        j120 &  old.math.nsc.ru/LBRT/k4/j120-hrs-Goncharov.txt \\ \hline
    \end{tabular}
\end{table}

Average processing time is 16 seconds for $\lambda=50000$ and 150
seconds for $\lambda=500000$ (for the instances with 120 activities).

\section{Conclusion}
Coelho and Vanhoucke \cite{CoelhoVanhoucke2020} encouraged researchers
 to focus their attention on developing radically new algorithms
 to solve RCPSP, rather than incrementally improving current
 algorithms that can solve existing instances of RCPSP only slightly better.
Such approaches await their researchers, and this work
is without a doubt part of the old conservative school of thought.
Nevertheless, we proposed a hybrid GA and NS algorithm for the RCPSP.
It incorporating two crossover techniques and two neighborhood variations.
We utilize a ranking system based on the importance of each resource.
Numerical experiments were carried out using datasets from the PSPLIB online library.
Computational experiments reveal that the suggested algorithm is a competitive heuristic,
outperforming multiple heuristics documented in the literature.
The best heuristic solutions were improved
for some instances from the j90 and j120 dataset.

%\clearpage

\section*{Acknowledgement}
The study was carried out within the framework of the state contract of the Sobolev
Institute of Mathematics (project FWNF-2022-0019)

%
% ---- Bibliography ----
%
% BibTeX users should specify bibliography style 'splncs04'.
% References will then be sorted and formatted in the correct style.
%
% \bibliographystyle{splncs04}
% \bibliography{mybibliography}

\begin{thebibliography}{48}
\bibitem{Abdolshah2014}
Abdolshah, M.:
A review of resource-constrained project scheduling problems
(RCPSP) approaches and solutions,
International Transaction Journal of Engineering, Management, \& Applied Sciences \& Technologies,
\textbf{5}(4), 253--286 (2014).

\bibitem{Alcaraz2006}
Alcaraz, J.,  Maroto, C.:
A hybrid genetic algorithm based on intelligent encoding for project scheduling,
In Jyzefowska J., \&  Weglarz J. (Eds.), Perspectives in modern project scheduling,
Boston: Springer, 249--274 (2006).

\bibitem{Berthaut2018}
Berthaut, F., Pellerin, R., Hajji, A., Perrier, N.:
A path relinking-based scatter search for the resource-constrained project scheduling problem,
International Journal of Project Organisation and Management, \textbf{10}(1), 1--36 (2018).

\bibitem{Blazewicz}
Bla\.{z}ewicz,~J., Lenstra,~J.K., Rinnoy Kan,~A.H.G.:
Scheduling Subject to Resource Constraints: Classification and
Complexity, Discrete Applied Math. \textbf{5}(1), ~11--24 (1983).

\bibitem{Brucker1999}
Brucker,~P., Drexl,~A., M\"{o}hring,~R., at al.:
Resource-Constrained Project Scheduling: Notation, Classification,
Models, and Methods, Eur. J. Oper. Res. \textbf{112}(1), ~3--41 (1999).

\bibitem{Carlier2009}
Carlier, J., Moukrim, A., Xu, H.:
A memetic algorithm for the resource constrained project scheduling
problem, In Proc. of the international conference on industrial
 engineering and systems management, IESM (2009).

\bibitem{ChenLiu2014}
Chen, D., Liu, S., Qin, S.:
Memetic algorithm for the resource-constrained project scheduling problem,
 In Proceeding of the 11 th world congress on intelligent control and automation,
WCICA,  4991--4996. IEEE (2014).

\bibitem {CoelhoVanhoucke2020}
Coelho, J., Vanhoucke, M.: Going to the core
of hard resource-constrained project scheduling instances.
Computers \& Operations Research, 121, Article 104976. (2020).

\bibitem{Czogalla2009}
Czogalla, J., Fink, A.:
Particle swarm topologies for resource constrained project scheduling,
In Krasnogor N. et al. (Eds.), Nature inspired cooperative strategies for optimization, 61--73.
Berlin Heidelberg: Springer-Verlag (2009).

\bibitem{Debels2006}
Debels, D., De Reyck Leus, B.R., Vanhoucke, M.:
A Hybrid Scatter  Search Electromagnetism Meta-Heuristic
for Project Scheduling,
Eur. J. Oper. Res.  \textbf{169}, 638--653 (2006).

\bibitem{DebelsVanhoucke2007}
Debels, D., Vanhoucke, M.:
Decomposition-based Genetic Algorithm
for the Resource-Consrtained Project Scheduling Problem,
Oper. Res.  \textbf{55}, ~457--469 (2007).

\bibitem{Elsayed2017}
Elsayed, S., Sarker, R., Ray, T., Coello, C. C.:
Consolidated optimization algorithm for resource-constrained
project scheduling problems, Information Sciences,  418--419, ~346--362 (2017).

\bibitem{FangWang2012}
Fang, C., Wang, L.:
An effective shuffled frog-leaping algorithm for resource-constrained project scheduling
problem, Computers \& Operations Research, \textbf{39}(5), 890--901 (2012).

\bibitem{Gagnon2004}
Gagnon, M., Boctor, F.F., d'Avignon, G.:
 A Tabu Search Algorithm for the Resource-constrained Project Scheduling Problem,
ASAC (2004).

\bibitem {Gim}
Gimadi, E.Kh.:
 On Some Mathematical Models and Methods for Planning Large-Scale Projects.
Models and Optimization Methods,  in Proc. AN USSR Sib. Branch, Math. Inst.,
Novosibirsk. Nauka \textbf{10}, 89--115 (1988).

\bibitem {GimSev2003}
  Gimadi E.Kh., Sevastianov S.V. On solvability of the project scheduling problem with accumulative resources of an
arbitrary sign. Operations Research Proceedings. Springer Verlag Berlin/Heidelberg/New York, 241--246 (2003).



\bibitem{Goncalves2011}
Goncalves, J., Resende, M.G.C, Mendes, J.:
 A Biased Random Key Genetic Algorithm with
Forward-Backward Improvement for
Resource-Constrained Project Scheduling Problem, J. Heuristics.
 \textbf{17}, 467--486 (2011).

\bibitem{Gon2017}
Goncharov, E.~N., Leonov, V.~V.:
 Genetic Algorithm for the Resource-Constrained Project Scheduling Problem,
Automation and Remote Control \textbf{78}(6) 1101--1114, (2017).

\bibitem {Gon2022a}
 Goncharov, E.~N. :
 An Improved Genetic Algorithm for the Resource-Constrained
 Project Scheduling Problem. In: Olenev, N., Evtushenko, Y.,
 Jacimovic, M., Khachay, M., Malkova, V., Pospelov, I. (eds)
 Advances in Optimization and Applications. OPTIMA 2022.
 Communications in Computer and Information Science, vol 1739.
 Springer, Cham. (2022).
 https://doi.org/10.1007/978-3-031-22990-9\_3

\bibitem {Gon2022}
Goncharov, E.N.
A Local Search Algorithm for the Resource-Constrained Project
Scheduling Problem.  Journal of Applied and Industrial Mathematics.
16, 672--683 (2022).
https://doi.org/10.1134/S1990478922040081

\bibitem{HarBris2010}
Hartmann, S., Briskorn, D.:
  A Survey of Variants and Extentions
of the Resource-Constrained Project Scheduling Problem,
 Eur. J. Oper. Res. \textbf{207}, ~1--14 (2010).

\bibitem{HarBris2021}
Hartmann, S., Briskorn, D.: An updated survey of variants and
extensions of the resource-constrained project scheduling problem.
European Journal of 297(1), 1--14. (2021).

\bibitem {Herroelen1998a}
Herroelen, W., Demeulemeester, E., De Reyck, B.:
 A Classification Scheme for Project Scheduling, Weglarz J. (Ed.).
Project Scheduling-Recent Models, Algorithms and
Applications, International Series in Operations Research and
Management Science.  Kluwer Acad. Publish.,
Dordrecht: \textbf{14}(1), 77--106 (1998).

\bibitem {Ismail2012}
Ismail, I. Y., Barghash, M. A.:
 Diversity guided genetic algorithm to solve the resource constrained project
scheduling problem,
International Journal of Planning and Scheduling, \textbf{1}(3), 147--170 (2012).

\bibitem {Jedrzejowicz2006}
Jedrzejowicz, P., Ratajczak, E.:
 Population learning algorithm for the resource-constrained project
scheduling, In Jyzefowska J., \& Weglarz J. (Eds.),
Perspectives in modern project scheduling, 275--296.
Boston: Springer (2006).

\bibitem{KolHar2006}
Kolisch, R., Hartmann, S.:
 Experimental Investigation of
Heuristics for Resource-Constrained Project Scheduling: An Update,
Eur. J. Oper. Res. \textbf{174}, ~23--37 (2006).

\bibitem{KolHar1999}
Kolisch, R., Hartmann, S.:
Heuristic Algorithms for Solving the Resource-Constrained Project Scheduling
Problem: Classification and Computational Analysis,
 Weglarz J., (ed). Project scheduling: Recent models, algorithms and applications.
 Kluwer Acad. Publish.,~147--178 (1999).

\bibitem {PSPLIB}
Kolisch, R., Sprecher, A.:
  PSPLIB -- a Project Scheduling Problem Library,
Eur. J. Oper. Res. \textbf{96}, 205--216 (1996).
(downloadable from http://www.om-db.wi.tum.de/psplib/)

\bibitem {Koulinas2014}
Koulinas, G., Kotsikas, L., Anagnostopoulos, K.:
 A particle swarm optimization based hyper-heuristic
algorithm for the classic resource constrained project scheduling problem,
Information Sciences, \textbf{277} 680--693 (2014).

\bibitem {Lim2013}
Lim, A., Ma, H., Rodrigues, B., Tan, S.T., Xiao, F.:
  New meta-heuristics for the resource-constrained project
scheduling problem,
Flexible Services and Manufacturing Journal, \textbf{25}(1-2), 48--73 (2013).

\bibitem{Mendes2009}
Mendes, J.J.M., Goncalves, J.F., Resende, M.G.C.:
  A Random Key Based Genetic Algorithm for the
Resource Constrained Project Scheduling Problem,
Comput. Oper. Res.  \textbf{36}, 92--109 (2009).

\bibitem {Mobini2009}
Mobini, M.D.M., Rabbani, M., Amalnik, M.S., at al.:
 Using an Enhanced Scatter Search Algorithm for a
Resource-Constrained Project Scheduling Problem, Soft Computing.
 \textbf{13} 597--610 (2009).

\bibitem {Mobini2011}
Mobini, M., Mobini, Z., Rabbani, M.:
  An artificial immune algorithm for the project scheduling problem
under resource constraints, Applied Soft Computing,
\textbf{11}(2), 1975--1982 (2011).

\bibitem{Palpant2004}
Palpant, M., Artigues, C., and Michelon, P.:
 Solving the resource-constrained
project scheduling problem with large neighborhood search,
Ann Oper Res \textbf{131} 237--257 (2004).

\bibitem{Pellerin2020}
Pellerin, R., Perrier, N., Berthaut, F.,:
 LSSPER: A survey of hybrid metaheuristics for the resource-constrained
project scheduling problem,
Eur. J. Oper. Res.  \textbf{280}, 2, 395--416 (2020).

\bibitem{ProonJin2011}
 Proon, S., Jin, M.:
 A Genetic Algorithm with Neighborhood Search
for the Resource-Consrtained Project Scheduling Problem,
Naval Res. Logist. \textbf{58} , 73--82 (2011).

\bibitem{Ranjbar2009}
Ranjbar, M., Kianfar, F.:
  A hybrid scatter search for the RCPSP,
Scientia Iranica, \textbf{16}(1), 11--18 (2009).

\bibitem{Paraskevopoulos2012}
Paraskevopoulos, D. C., Tarantilis, C. D., Ioannou, G.:
 Solving project scheduling problems with resource
constraints via an event list-based evolutionary algorithm,
Expert Systems with Applications, \textbf{39}(4), 3983--3994 (2012).

\bibitem{Valls2005}
Valls, V., Ballestin, F., Quintanilla, M.S.:
  Justification and RCPSP: a Technique that Pays,
Eur. J. Oper. Res. \textbf{165}, 375--386 (2005).

\bibitem{Valls2008}
Valls, V., Ballestin, F., Quintanilla, S.:
A Hybrid Genetic Algorithm for the Resource-Consrtained
Project Scheduling Problem,
Eur. J. Oper. Res.  \textbf{185}(2), 495--508 (2008).

\bibitem {Valls2004}
Valls, V., Ballestin, F., Quintanilla, S.:
  A Population-based Approach
to the Resource-Constrained Project Scheduling Problem,
Annals of Operations Research \textbf{131}, 305--324 (2004).

\bibitem {Vanhoucke2012}
Vanhoucke, M.:
  Resource-constrained project scheduling,
In Project Management with Dynamic Scheduling, 107--137,
Berlin, Heidelberg: Springer-Verlag. (2012)

\bibitem {VanhouckeCoelho2024}
Vanhoucke, M. and Coelho, J., 2024, A matheuristic for the
resource-constrained project scheduling
problem, European Journal of Operational Research,
319(3), 711--725 (doi: 10.1016/j.ejor.2024.07.016)


\bibitem{WangChen2010}
Wang Chen, Yan-jun Shi, Hong-fei Teng, at al.:
  An Efficient Hybrid Algorithm for Resource-Constrained Project
Scheduling,  Inf. Sci.  \textbf{180}(6), 1031--1039 (2010).

\bibitem{WangLi2010}
Wang, H., Li, T., Lin, T.:
  Efficient genetic algorithm for resource-constrained project scheduling problem,
Transactions of Tianjin University, \textbf{16}(5), 376--382 (2010).

\bibitem {Weglarz1999}
Weglarz, J.:
  Project scheduling. Recent models, algorithms and applications,
Boston: Kluwer Acad. Publ. (1999).

\bibitem {YingLinLee2009}
Ying, K. C., Lin, S. W., Lee, Z. J.:
  Hybrid-directional planning: improving improvement heuristics for
scheduling resource-constrained projects,
International Journal of Advanced Manufacturing Technology, \textbf{41}(3--4),
 358--366 (2009).

\bibitem {Zamani2013}
Zamani, R.:
  A competitive magnet-based genetic algorithm for solving the resource-constrained project
scheduling problem,
European Journal of Operational Research, \textbf{229}(2), 552--559 (2013).

\bibitem {ZhengWang2015}
Zheng, X., Wang, L.:
  A multi-agent optimization algorithm for resource constrained project scheduling
problem, Expert Systems with Applications, \textbf{42}(15--16), 6039--6049 (2015).

\end{thebibliography}
%

\end{document}